\documentclass[12pt]{amsart}
\usepackage{amsthm}
\usepackage{amsmath}
\usepackage{amsfonts}

\title[Planar configurations and gkz-rationality]{\bf
Planar configurations of lattice
vectors and gkz-rational toric fourfolds in ${\mathbb
P}^6$}

\author{Eduardo Cattani and Alicia Dickenstein}
\address{Eduardo Cattani: Department of Mathematics
and Statistics. University
of Massachusetts. Amherst, MA 01003, USA}
\email{cattani@math.umass.edu}
\thanks{E. Cattani is partially supported by NSF Grant
DMS--0099707}

\address{Alicia Dickenstein: Departamento~de
Matematica, FCEyN.
Universidad de Buenos Aires. (1428) Buenos Aires,
Argentina}
\email{alidick@dm.uba.ar}
\thanks{A. Dickenstein is partially supported  by 
UBACYT X052 and ANPCYT 03-6568,
Argentina.}

\newcommand{\baseRing}[1]{\ensuremath{\mathbb{#1}}}
\newcommand{\Z}{\baseRing{Z}}
\newcommand{\R}{\baseRing{R}}
\newcommand{\C}{\baseRing{C}}
\newcommand{\N}{\baseRing{N}}
\newcommand{\Q}{\baseRing{Q}}

\newcommand{\CP}{\baseRing{P}}

\def\pd#1{ \partial_{#1} }

\theoremstyle{plain}
\newtheorem{theorem}{Theorem}[section]
\newtheorem{lemma}[theorem]{Lemma}
\newtheorem{corollary}[theorem]{Corollary}

\newtheorem{proposition}[theorem]{Proposition}

\newtheorem{assumption}[theorem]{Assumption}

\theoremstyle{definition}

\newtheorem{definition}[theorem]{Definition}
\newtheorem{remark}[theorem]{Remark}

\newtheorem{example}[theorem]{Example}

\numberwithin{equation}{section}

\newcommand{\Script}[1]{\ensuremath{{\mathcal{#1}}}}

\newcommand{\PP}{\Script{P}}

\newcommand{\EE}{\Script{E}}

\newcommand{\LL}{\Script{L}}

\begin{document}

\begin{abstract}

We introduce a notion of  balanced configurations of
vectors. This is motivated by the study of rational
$A$-hypergeo\-me\-tric functions in
the sense of Gelfand, Kapranov and Zelevinsky.
We classify  balanced configurations of seven
plane vectors up to $GL(2,\R)$-equivalence and
deduce  that the only gkz-rational toric four-folds
in $\CP^6$ are
those varieties associated with an essential Cayley
configuration. We show that in this case, all rational
$A$-hyper\-geo\-metric functions
may  be described in terms of toric residues.  This
follows from
studying a suitable  hyperplane arrangement.
\end{abstract}

\maketitle

\section{Introduction}
A configuration $\{b_1,\dots,b_n\}$ of vectors
in $\R^m$ is said to be {\em balanced} if for 
every index set $J = \{1\leq  j_1 < \cdots <  j_{m-1}\leq n\}$
the multiset 
$$\{\det(b_{j_1},\dots, b_{j_{m-1}},b_i); i\not\in J \}$$
is symmetric around the origin.  

Balanced configurations in the plane with at most
six vectors
have been classified in \cite{rhf}.  Here we consider
the case of seven planar  vectors and show in
Corollary~\ref{balancedconf} that
 there are only five balanced  configurations 
up to  $GL(2,\R)$-equivalence. When the configuration
is uniform, i.e. every pair of vectors
is linearly independent, it is equivalent to a regular
heptagon.
In fact, if an arbitrary uniform configuration  
of plane vectors is balanced, 
then it is $GL(2,\R)$-equivalent to a regular 
$(2k +1)$-gon. 
 This result was conjectured in the preprint 
 version of this
 paper and recently verified by N. Ressayre \cite{ressayre}.

Our interest in  balanced configurations stems from
their relationship with multivariable
hypergeometric functions in the sense of Gel'fand,
Kapranov, and Zelevinsky
\cite{gkz89, gkz90}. These functions are solutions of
a regular, holonomic system
of partial differential equations associated with a
configuration
$A\subset\Z^d$ and a homogeneity vector $\alpha \in
\C^d$.
They include, as particular examples, the classical
Gauss hypergeometric
function, as well as the multivariable generalizations
of Appell, Horn, and Lauricella~\cite{erdelyi}.

The combinatorics of the configuration $A$ plays a
central role in the study of the $A$-hypergeometric
system: from the construction of  series solutions 
associated with regular triangulations of the convex hull 
of $A$~\cite{gkz89}, to the recent
Gr\"obner deformation methods of Saito, Sturmfels and
Takayama~\cite{sst2}.  Conversely, the existence of
rational solutions
imposes strong combinatorial restrictions on $A$.

For appropriate
integer homogeneities,
every configuration $A$ admits polynomial solutions; 
indeed, they are
related to the
associated integer programming problem ~\cite{sst}. 
Similarly, for
every $A$ and suitable $\alpha$, there exist Laurent
polynomial
solutions and, for projective curves, their number
reflects algebraic properties
of the toric ideal associated with $A$~\cite{cdd}.  We
should also mention in this context, the
classical work of Hermann Schwarz on algebraic
solutions of Gauss'
hypergeometric equation~\cite[\S 10.3]{hille}

In joint work with Bernd Sturmfels \cite{rhf, binom},
we have considered
the combinatorial restrictions on $A$ imposed by the
existence of
$A$-hyper\-geo\-metric rational functions other than
Laurent polynomials.
It is well-known that the irreducible components of
the singular
locus of the $A$-hyper\-geo\-metric system are defined
by {\em facial discriminants}, i.e. the sparse discriminant 
of a subconfiguration given
by the intersection of $A$ with a face of its convex hull
\cite{gkz89,gkzbook}.
In particular, such discriminants are the possible factors of the
denominator of
a rational $A$-hypergeometric function.  Following
\cite{rhf} we say
that $A$ is {\em gkz-rational} if the discriminant
$D_A$ is not
a monomial and there exists a (non-zero) rational
$A$-hypergeometric function $f$ whose denominator is a
multiple of  $D_A$. It is conjectured in
\cite[Conjecture~1.3]{rhf} that a configuration is
gkz-rational if and only if it is affinely equivalent to an {\it
essential Cayley
configuration\/} (c.f. (\ref{Cayley})).  This conjecture
has been verified in \cite{rhf} when the projective
toric variety
associated with  $A$ is a hypersurface or has
dimension at most  three.  Indeed, in the hypersurface
case there is a direct correspondence between 
gkz-rationality and balanced
configurations.  Let  $A$ be a $d \times (d+1)$-matrix
of rank $d$ and
$(b_1,\dots,b_{d+1})$  a $\Z$-generator of the
integral kernel of $A$. Then \cite[Theorem~2.3]{rhf} 
asserts that $A$ is gkz-rational if and
only if the multiset $\{b_1,\dots,b_{d+1}\}$ is
symmetric around the origin.

In this paper we study the first open case, that
of codimension-two toric subvarieties of $\CP^6$, in order
to reveal the combinatorial and analytic problems that arise. 
Our starting point are the results in \cite{elim} about the sparse
discriminant of codimension-two configurations.  They allow us to
show that in a gkz-configuration every {\em non-splitting 
circuit} (cf. Definition \ref{defessentialcircuit})
must be balanced.  Via Gale duality we are
then led  to the classification of balanced
configurations of seven lattice vectors in the plane.

In \S 3 we deduce from the classification
Theorem~\ref{classification}
  that Conjecture~1.3 in \cite{rhf} holds for toric
four-folds in $\CP^6$.
That leads naturally to the question of describing the
space of
$A$-hypergeometric rational functions associated with
the corresponding
essential Cayley configuration.  We show that, as
predicted by
\cite[Conjecture~5.7]{rhf}, a suitable derivative of
such a function
can be realized as a toric residue in the sense
of~\cite{compo, cox}.
As in \cite{binom}, this is done by studying  an
appropriate oriented hyperplane arrangement.
We conclude the paper with an example illustrating
these constructions
for the configuration associated with Appell's
classical function $F_2$.

\smallskip

\noindent{\bf Acknowledgements:}  We are grateful to
Bernd Sturmfels for many valuable discussions and to the
two anonymous referees for their comments and suggestions
for improvement.  This work was
completed while Alicia Dickenstein was visiting the University of
Stockholm with the support of the Swedish Science Council.  She is
grateful to the Department of Mathematics and to Mikael Passare for
their hospitality.

\section{Classification of two dimensional configurations}
\label{geom}

Let $A =  \{a_1,\dots,a_{n}\} \subset
\Z^d$, be an integral configuration spanning $\Q^d$.
We denote by $m:=n-d$ the {\em codimension } of $A$.
Let $B\in \Z^{n\times m}$ be a matrix, well defined up
to right multiplication by elements of $GL(m,\Z)$, whose
columns $\nu_1,\dots, \nu_{m}$ are a
$\Z$-basis of the lattice
\begin{equation}\label{lattice}
\LL \ :=  \ \{v\in \Z^n : A\cdot v =   0\}.
\end{equation}
Let $b_i =   (b_{i1},\dots,b_{im}) \in \Z^m$ be the
$i$-th row of $B$.
We shall also denote by $B$ the configuration
$\{b_1,\dots,b_n\}\in \Z^m$ and refer to it as a
{\em Gale dual} of $A$.

Note that the vector $(1,\dots,1)$ is in the row span
of $A$ if and only if
\begin{equation}\label{regularity}
\sum_{i=  1}^n b_i \ =  \ 0.
\end{equation}
Configurations satisfying (\ref{regularity}) are
called {\em nonconfluent\/}.  This terminology will be
explained in \S 3 when we discuss the hypergeometric
system associated with a configuration.

Since in this paper we are ultimately interested in 
studying rational bivariate
$A$-hypergeometric functions we shall assume
from now on that $A$ is a nonconfluent codimension-two
configuration.
We  also make, without loss
of generality, two simplifying assumptions:

\begin{assumption}\label{pyramid}
No hyperplane contains all but one
point of $A$. Equivalently, all $b_i
\not=  0.$
\end{assumption}

When this assumption is not satisfied, the
discriminant
$D_A=  1$, and so $A$ cannot be gkz-rational.

\begin{assumption}\label{distinctpoints}
The configuration
$A$ consists of distinct points.
Du/-al/-ly,  no line contains all but
two of the vectors in  $B$.
In particular, we may assume from now on that $n\geq
4$.
\end{assumption}

In fact, the study of $A$-hypergeometric functions
associated to configurations with repeated points
reduces trivially to the study of hypergeometric
functions associated to the configuration of the distinct 
points in $A$.

\vskip .3cm
Given an index set $I \subset \{1,\dots,n\}$, let
$A(I)$ denote the 
subset $A(I) = \{a_i\in A; i\in I\}$.  
Recall that a {\em circuit} in $A$ is a
minimally dependent subset of
$A$. A circuit $A(I)$ defines a codimension-one
configuration so its
Gale dual is given by a multiset
$  \{c_i \ , \ i \in I \}\subset \Z\backslash \{0\}$,
well defined up to multiplication by $-1$. These scalars
define a $\Z$-minimal relation
$$\sum_{i\in I} c_i a_i =   0\,.$$
The vector $c \in \Z^n$ whose $i$-th component is
$c_i$ and whose $j$-th component vanishes for all
$j\in I^c$, belongs to $\LL$ and, consequently, may
be written as $c =   \gamma_1 \nu_1 + \gamma_2 \nu_2$.
This means
that the vectors $\{b_j, \  j\in I^c\}$  lie
in the line orthogonal to
$(\gamma_1,\gamma_2)\in \Z^2$.   A subset $B(J) = \{b_j; j\in J\}$
is a {\em cocircuit}   if and only if it consists of all
vectors of $B$  in the same
line through the origin.  Clearly, $B(J)$ is a cocircuit
if and only if
$A(J^c)$ is a circuit. In this case, given $j\in J$, 
the multiset 
\begin{equation}\label{galecircuit}
\{\det(b_i,b_j)\,;\ i\notin J \}
\end{equation}
is well defined up to multiplication by a non-zero constant and
agrees, again up to constant, with the Gale dual  of
the circuit $A(J^c)$.  In particular, the following definition
is independent of the choice of $j\in J$.

\begin{definition}
A cocircuit $B(J)$ is said to be {\em balanced}
if for some (all) $j\in J$, the multiset
$\{\det(b_i,b_j)\,;\ i\notin J \}$ 
is symmetric around the origin.  $B$ is said to
be {\em balanced}
if all its cocircuits are balanced.  Also, a circuit
$A(I)$ is balanced if $B(I^c)$ is so or, equivalently, if
its Gale dual is symmetric around the origin. 
\end{definition}

A fundamental property of Gale duality 
(cf. \cite[Proposition 9.1.5]{matroids})
states that  the convex hull of a subset
$A(I)$ is a face of the convex hull of
$A$ if and only if the origin lies in the relative
interior of ${\rm conv}(B( I^c))$.  Hence, in the
codimension
two case, a circuit $A(I)$ is contained in a
face of ${\rm conv}(A)$ if and only if the cocircuit
$B(I^c)$ contains
 at least two vectors of $B$ in opposite rays.  
Moreover, if $B(J)$ is a cocircuit and
\begin{equation}\label{essentialcircuit}
\sum_{j\in {J}} b_j =  0,
\end{equation}
the configuration $A$ can be written as $A = A(J) \cup
A(J^c)$ with $A(J)$ and
$A(J^c)$ lying in parallel flats. 

\begin{definition} \label{defessentialcircuit}
A cocircuit $B(J)$ is 
{\em splitting} if  equality (\ref{essentialcircuit}) holds.
Otherwise we  call $B(J)$ {\em non-splitting}.  We apply
the same terminology to the corresponding circuit $A(J^{c})$.
\end{definition} 

As Lemma~\ref{discspec} shows, non-splitting circuits
behave well
relative to discriminants and specialization of
variables.  Moreover,
Theorem \ref{generic} implies that in order to
classify gkz-rationality, we must study  configurations
all of whose non-splitting circuits are balanced.

\smallskip

We give now two definitions that will help us organize
the arguments. Let $B\subset \R^2$ be a finite vector
configuration.

\begin{definition}\label{unidef}
We say that $B$ is
{\em uniform} if for every pair $i\not=   j$, the
vectors $b_i, b_j$ are linearly independent.
\end{definition}

Note that in a uniform configuration every cocircuit
is  non-splitting. We also observe that any
cocircuit whose complement consists of an odd number of
vectors, is automatically unbalanced.  In particular,
a uniform balanced configuration must necessarily consist
of an odd number of points.

\begin{definition}\label{irredef}
We say that $B$ is {\em irreducible} if it is
not possible to write it as a disjoint union $B =  
B_1 \cup B_2$,
where $B_1$, $B_2$ are proper nonconfluent
subconfigurations.
\end{definition}

We make the following convention: each time we write
$B$ as a union
of subconfigurations, we assume that we are in the
conditions of
the previous definition, i.e., that the union is
disjoint and the
subconfigurations are proper and with zero sum. In
particular, each
subconfiguration contains at least two vectors.

\begin{remark} \label{necess}  If a
configuration
satisfies Assumption \ref{pyramid}, all
subconfigurations do as well.
Assumption \ref{distinctpoints}, on the contrary, is
not  hereditary.
\end{remark}

\smallskip

For  $4 \leq n \leq 6$,
the geometric classification of two dimensional Gale 
configurations satisfying
Assumptions \ref{pyramid} and \ref{distinctpoints}
follows from the results in \cite{rhf} which we now recall.

\begin{itemize}

\smallskip

\item If $n=  4$, every configuration $B$ is
irreducible and uniform
and therefore  every cocircuit is
non-splitting and unbalanced.

\smallskip

\item For $n=  5$,  \cite[Lemma~4.2]{rhf} implies that
either

\begin{itemize}
\item[a)]$B$ is irreducible and contains a
non-splitting unbalanced cocircuit, or
\item[b)]
$B$ is irreducible and all its cocircuits are
balanced, in which case
it is $GL(2,\R)$-equivalent to the regular
pentagon, or
\item[c)] $B$ is reducible; hence it is the union of a
two-vector
configuration $B_1$ and a three-vector configuration
$B_2$ not contained in a line. If $B$ is
integral  then it  is Gale dual to an essential Cayley
configuration
(\ref{Cayley}).
\end{itemize}

\smallskip

\item Suppose $n=  6$. Then we have the
following possibilities:
\begin{itemize}
\item[a)] $B$ contains a non-splitting unbalanced
cocircuit, or
 \item[b)] $B$ is $GL(2,\R)$-equivalent to $B_1 \cup
B_2$,
where 
$$B_i= \{ \lambda_i e_1, \lambda_i e_2, - \lambda_i
e_1 - \lambda_i e_2\},$$  
and $\lambda_1 +
\lambda_2 \not= 0.$ In this case,  all cocircuits are
balanced.
\item[c)] $B$ is a reducible configuration of the form
$B =   B_1\cup B_2 \cup B_3$ or $B=  B_1 \cup B_2$
and each $B_i$ is a cocircuit.
\end{itemize}
\end{itemize}

The following is the main theorem in this section.

\begin{theorem}\label{classification}
Let $B =   \{b_1,\dots,b_7\}\subset \R^2 $ be a
nonconfluent configuration of
seven non-zero vectors in the plane.
Assume moreover that no line contains all but two of
the vectors.
Then, exactly one of the following holds:
\begin{itemize}
\item[a)] $B$ contains a non-splitting unbalanced
cocircuit.
\item[b)] $B$ is irreducible and all its cocircuits
are balanced,
in which case it is
$GL(2,\R)$-equivalent to a regular heptagon.
\item[c)] $B$ is a reducible configuration of the form
$B =   B_1\cup B_2 \cup B_3$ or $B=  B_1 \cup B_2$
with each $B_i$ a cocircuit.
\item[d)] $B$ is a reducible configuration of the form
$B =   B_1\cup B_2 \cup B_3$, with one  $B_i$ having
three vectors which generate $\R^2$.
\item[e)] $B$ is a reducible configuration of the form
$B =   B_1\cup B_2$, where $B_1$ is a pair of opposite
vectors and $B_2$ is $GL(2,\R)$-equivalent to a
regular pentagon.
\end{itemize}
\end{theorem}

Recall that an arbitrary configuration $B$ of planar
vectors
is called balanced if all its cocircuits are balanced.

\begin{corollary}\label{balancedconf}
Let $B$ be a  balanced configuration of 
seven non-zero vectors in the plane.  Then $B$ is 
$GL(2,\R)$-equivalent to one of the following:
\begin{itemize}
\item[i)] A regular heptagon,
\item[ii)] a configuration $B = B_1 \cup B_2$, where
$B_2$ is  a regular pentagon, and $B_1 = \{\lambda
c,-\lambda c\}$,
$c\in B_2$, $\lambda\in \R^*$,
\item[iii)] the configuration $B = \{e_1,e_2,
-e_1-e_2,\lambda e_1, -\lambda e_1,
\mu e_2, -\mu e_2\}$, where $\lambda, \mu$ are
non-zero scalars,
\item[iv)] the configuration $B = \{e_1,e_2,
-e_1-e_2,\lambda e_1, -\lambda e_1,
\mu e_1, -\mu e_1\}$, where $\lambda, \mu$ are
non-zero scalars,
\item[v)] a configuration of seven vectors in one
line.
\end{itemize}
\end{corollary}

\begin{proof} A configuration $B$ all of whose vectors
are contained
in a line is vacuously balanced.  On the other hand,
if $B$ spans
the plane, then a balanced configuration is
necessarily nonconfluent.
If, in addition, $B$ satisfies Assumption~\ref{distinctpoints}
 then it must be equivalent to one of the configurations a)
-- e) in Theorem~\ref{classification}.  Of those, the
configurations of
type a) and c) cannot be balanced.  If a configuration
of type e) is balanced then it must be as in ii) since
otherwise  $B_1$ would be an unbalanced cocircuit.  

Suppose now that $B$ is as in d).  Then, modulo
$GL(2,\R)$, we may
assume that $B_1 = \{e_1,-e_1\}$, $B_2 =
\{e_2,-e_2\}$.  Since
the corresponding cocircuits must be balanced, this
implies that
$B_3 = \{\lambda e_1, \mu e_2, -\lambda e_1 - \mu
e_2\}$ which, 
after rescaling, yields the configuration iii).

It remains to consider the case when a line contains
at least
five vectors of $B$.  Clearly, since $B$ is
nonconfluent, we must
consider
the case when all seven vectors lie in a line (case
v), or there
are exactly five vectors in one line.  If the other
two vectors
were also collinear then the corresponding cocircuit
would be
unbalanced.  So we may assume that the remaining two
vectors are
linearly independent and that $B$ is
$GL(2,\R)$-equivalent to
$B = \{\lambda_1 e_1, \lambda_2 e_1, \lambda_3 e_1,
\lambda_4 e_1, \lambda_5 e_1,
\lambda_6 e_1 + \gamma e_2, e_2\}, \lambda_6\not=0$. 
Condition \ref{regularity} implies $\gamma =
-1$, while the fact that the cocircuit defined by
$e_2$ is balanced implies that,
after reordering if necessary, $\lambda_1 =
-\lambda_2$, $\lambda_3 = -\lambda_4$,
and $\lambda_5 = -\lambda_6$.  Rescaling yields
configuration iv).
\end{proof}

The proof of Theorem~\ref{classification} will follow
from a series of Lemmas. In their statements, a)--e)
refer to the cases described in the theorem. 
Throughout, all
 configurations are supposed to verify Assumptions
\ref{pyramid} and \ref{distinctpoints}.
Our first goal is to show that if $B$ in not uniform
then either it
contains a non-splitting unbalanced cocircuit or is as
in  c), d) or e).

\begin{lemma}\label{reduction}
Suppose that $B$ is a reducible configuration
such that $B =   B_1 \cup B_2$, where
$B_1 =   \{b,-b\}$. Then $B$ is as in a), c), d) or e).
\end{lemma}

\begin{proof}
By Assumption~\ref{distinctpoints}, $B_2$ is a
$5$-vector configuration
not contained in a line. Hence, either
\begin{itemize}
\item $B_2$ contains a non-splitting unbalanced
cocircuit  associated with
a line $L$. Then $L$ also defines a non-splitting
unbalanced cocircuit of the total configuration
$B$, or
\item $B_2$ is a reducible configuration
which does not satisfy Assumption~\ref{distinctpoints}
and  we are in case c), or
\item $B_2$ is a reducible configuration
satisfying Assumption~\ref{distinctpoints}.
Then so is $B$ and we are in case  d), or
\item $B_2$ is equivalent to a regular pentagon and
$B$ is as in e).
\end{itemize}~\end{proof}

\begin{lemma}\label{evennumber}
Suppose that there is a line $L$ containing an even
number of vectors of $B$.
Then $B$ is as in a), c), d), or e).
\end{lemma}

\begin{proof}
Since $L$ contains an even number of vectors, the
associated cocircuit
is unbalanced. If it is non-splitting, then we are in
case a). Otherwise,
 the vectors in $B\cap L$ must add up to zero. If $L$
contains
only two vectors then we are in the case covered by
Lemma~\ref{reduction}. Hence we may assume that $L$
contains
$4$ vectors whose sum is zero. If the remaining $3$
vectors lie on a line then we are in the case described by c);
otherwise, there is a line
$L'$ containing a single vector of $B$ and,
consequently, the associated
cocircuit is non-splitting. If it is unbalanced we are
in case a) while if it
is balanced then necessarily there exist $b_1, b_2\in
L$ such that
$b_1 =   -b_2$ and we are again in the case covered by
Lemma~\ref{reduction}.
\end{proof}

\smallskip

We can now restrict ourselves to configurations $B$
where each line
contains an odd number of vectors of $B$. This means,
in particular,
that $B$ is distributed in an odd number of lines.

\begin{lemma}\label{threelines}
Suppose that the configuration $B$ is distributed in
$3$ lines, each
containing an odd number of vectors. Then
$B$ contains a  non-splitting unbalanced cocircuit,
or it corresponds to cases c) or d).
\end{lemma}

\begin{proof} By Assumption~\ref{distinctpoints}, no
line may contain $5$ vectors of $B$.
Hence we need only consider
the case when one of the lines, say $L_1$, contains only
one vector, while
the other two $L_2, L_3$ contain three vectors each.
Let $Z_i\subset B$ denote
the cocircuit associated with $L_i$. Because of
Assumption~\ref{pyramid}, one of the cocircuits $Z_2$,
$Z_3$ must be
non-splitting. Assume $Z_2$ is non-splitting. If it is
balanced,
then there exist elements $b_1, b_2 \in Z_3$ such that
$b_1 =  - b_2$; hence we are in the case of
Lemma~\ref{reduction} and,
since $B$ is contained in three lines, it is as in a), c) or d).
\end{proof}

\begin{lemma}\label{fivelines}
Suppose that the configuration $B$ is distributed in
$5$ lines, each containing an odd number of vectors. Then
either $B$ contains a non-splitting unbalanced cocircuit, or
is as in d), or as in e) and the line containing
$B_1$ contains one of the vertices of the pentagon.
\end{lemma}

\begin{proof} It suffices to show that either $B$
contains a
non-splitting unbalanced cocircuit or it reduces to
the case in
Lemma~\ref{reduction}. Since one of the lines must
contain three
vectors of $B$, we may assume that $b_1 =   e_1$, $b_2
=  x e_1$, $b_3 = ye_1$
with $x,y \not=   -1$ and $x \not=   -y$. We may also
assume that
$b_4 =   e_2$. Hence, if the cocircuit associated with
the line $\R\cdot e_2$ is
balanced, the configuration $B$ must be of the form:
\begin{equation*}\label{config5}
B^t\ =   \ \left(
\begin{array}{ccccccc}
1 & x & y & 0 & -1 & -x & -y\\
0 & 0 & 0 & 1 & * & * & * \\
\end{array}
\right).
\end{equation*}
If the cocircuit defined by $\R\cdot e_1$ is also
balanced then either
\begin{equation}\label{config5a}
B^t\ =   \ \left(
\begin{array}{ccccccc}
1 & x & y & 0 & -1 & -x & -y\\
0 & 0 & 0 & 1 & -1 & u & -u \\
\end{array}
\right),\ u\not=  0,\ \hbox{or}
\end{equation}
\begin{equation}\label{config5b}
B^t\ =   \ \left(
\begin{array}{ccccccc}
1 & x & y & 0 & -1 & -x & -y\\
0 & 0 & 0 & 1 & u & -u & -1 \\
\end{array}
\right),\ u\not=  0.
\end{equation}
Consider the configuration (\ref{config5a}). If the
cocircuit associated
with $\R\cdot (-1,-1)$ is balanced, we must
have $u =   y - x$.
But then, the cocircuit defined by the line $\R\cdot (-x,u)
= \R\cdot (-x,y-x)$ cannot be balanced since the multiset
\begin{equation*}  
\ \{y-x,\ x(y-x),\ y(y-x),\ x,\ -y,\
-(y-x)(x+y)\}
\end{equation*}
can never be symmetric around the origin given that
 $x,y \not=   0, -1$ since
the configuration satisfies Assumption~\ref{pyramid}.

It remains to consider the case described by
(\ref{config5b}).
It is easy to check that under our assumptions, if the
cocircuit associated with the line $\R\cdot(-y,-u)$ is
balanced then we must have $x =   1 + uy$. Consider then 
the cocircuit defined by $\R\cdot (-1,u)$. Checking again 
case-by-case we deduce
that if it is balanced then $u =   1 + uy$ and hence
$u=  x$.
We can then check that the last cocircuit, the one
defined by
$\R\cdot(-x,u)$, will be unbalanced unless $y =   -x$.
\end{proof}

Lemma~\ref{fivelines} completes the proof of
Theorem~\ref{classification}
for non-uniform configurations.  
The following result,
was
 conjectured in the preprint version of this paper
 and proved for the case of seven plane vectors using
 Gr\"obner basis computations.  
 Recently, 
N.~Ressayre \cite{ressayre} gave a proof for the general
case.  
For the sake of completeness,
 we include an adaptation of Ressayre's argument to our
 situation.
 
\begin{theorem}\label{uniform}
Let $B$ be a uniform balanced configuration in the plane. Then
$B$ is $GL(2,\R)$-equivalent to a (2k + 1)-gon.
\end{theorem}

\begin{proof}
We refer to \cite{ressayre} for the general case and
sketch an argument in the case of seven vectors.

Let $\PP$ denote the set of unordered pairs $\{(i,j): 1\leq i , j \leq
7\}$ and $$\PP^{k} \ =\ \{(i,j)\in \PP : \det(b_{k},b_{i}) =
-\det(b_{k},b_{j})\}.$$
Note that if $(i,j) \in \PP^{k}$ then $b_{k}$ lies in the line
spanned by $b_{i}+b_{j}$.  Hence, the uniformity assumption implies
that the sets $\PP^{k}$ are disjoint.  A cardinality argument shows
that 
$\PP = \displaystyle{\cup_{k=1}^{7} \PP^{k}}$.

We also note that since every cocircuit is balanced, for each
$k=1,\dots,7$, there exist three vectors $b_{i}$ such that
$\det(b_{i},b_{k})>0$ and three vectors with negative determinant.
Hence, each of the open hyperplanes determined by the line $\R\cdot b_{k}$
must contain exactly three vectors in the configuration.  
Suppose now that we index the vectors $b_{1},\dots,b_{7}$ 
counterclockwise, then
the vector $b_{k}$ must be contained in the open cone spanned by
$-b_{k-3}$ and $-b_{k+3}$, where we are indexing modulo $7$. This
implies that $b_{k}$ is the only vector in the configuration that
may lie in the line $\R\cdot(b_{k-3}+b_{k+3})$ and, consequently,
\begin{equation}\label{pp1}
    (k-3,k+3) \in \PP^{k}.
    \end{equation}
    
Similarly, we have 
\begin{equation}\label{pp2}
    (k-1,k+1) \in \PP^{k}.
    \end{equation}
Indeed, by symmetry it suffices to verify (\ref{pp2}) for the
case $k=2$.  Then $(1,3)$ must belong to either $\PP^{2}$,
since $b_{2}$ is in the open cone spanned by $b_{1}$ and $b_{3}$,
or $\PP^{5}$, or  $\PP^{6}$ since $b_{5}$ and $b_{6}$ are in the
cone spanned by $-b_{1}$ and $-b_{3}$.  However, it follows from
(\ref{pp1}) that $(1,2) \in \PP^{5}$ and $(2,3) \in \PP^{6}$.
Consequently, the last two cases are impossible and $(1,3)\in\PP^{2}$,
as claimed by (\ref{pp2}).  
By a process of
elimination we also have
\begin{equation}\label{pp3}
    (k-2,k+2) \in \PP^{k},
    \end{equation}
and the set $\PP^{k}$ is completely determined.

Modulo the action of $GL(2,\R)$ we may assume that 
$b_{1} = e_{1}$ and $b_{2} = e_{2}$.  Since
$\PP^{1} = \{(2,7), (3,6),(4,5)\}$ and 
$\PP^{2} = \{(1,3), (4,7),(5,6)\}$ the configuration
$B$ must be of the form

\begin{equation}\label{b1}
B^t\ =   \ \left(
\begin{array}{ccccccc}
1 & 0 & -1 & z & w & -w & -z\\
0 & 1 & x & y & -y & -x & -1 \\
\end{array}\right).
\end{equation}

Since $(1,2)\in \PP^{5}$ it follows that $b_{5}=\lambda(e_{1}+e_{2})$ 
and, therefore, $w=-y$.  Moreover,  $(3,7)\in \PP^{5}$
as well, so  that $x = -z$.  We also have that
 $(1,6)\in \PP^{7}$ which implies that $y = x^{2}-1$.
 Therefore
\begin{equation}\label{b2}
B^t\ =   \ \left(
\begin{array}{ccccccc}
1 & 0 & -1 & -x & 1-x^{2} & x^{2}-1& x\\
0 & 1 & x & x^{2}-1 & 1-x^{2} & -x & -1 \\
\end{array}\right).
\end{equation}
Finally, since $(1,7)\in \PP^{4}$ we have that 
$$ x^{3} + x^{2} - 2 x -1 =0.$$ 
This equation has three real roots: 
$x =   2 \cos(2k\pi/7)$, $k=  1,2,3$. A
straightforward argument shows
that the corresponding configurations are
$GL(2,\R)$-equivalent to a regular heptagon.
\end{proof}

\section{Classification of rational hypergeometric
functions}
\label{hyper}

We recall the definition of $A$-hypergeometric
functions
and refer to \cite{gkz89,gkz90,sst2} for their main
properties.

\begin{definition}\label{def:hypergeom}
Given an  integer $ d \times n$-matrix $A$  of
rank $d$ and  a  vector $\alpha \in \C^d$,
the {\it $A$-hypergeometric system} with parameter
$\alpha$  is  the left ideal
$H_A(\alpha)$ in the Weyl algebra
$ \C\langle x_1,\dots,x_n,\pd 1,\dots,\pd n\rangle$
 generated by the {\sl toric operators}
$\ \partial^u - \partial^v$, for all
$ u,v\in \N^n$ such that   $A\cdot u=  A\cdot v$,
and the {\sl Euler operators}
$ \,\sum_{j=  1}^n a_{ij} x_j \partial_j - \alpha_i
\,$ for $\, i=  1,\dots,d $.
A function $f(x_1,\dots,x_n)$, holomorphic in an open
set $U\subset \C^n$, is said to be {\it
$A$-hypergeometric of
degree $\alpha$} if it is annihilated by $H_A(\alpha)$.
\end{definition}

We say that the dimension of $A$ is $d-1$, i.e. the
dimension of the affine span of its columns, and $m = n -d$ is
its codimension.  When $m=1$, the study of the 
hypergeometric system may be reduced to the study of
a hypergeometric ordinary differential equation of
degree $d$.  The singularities of this equation are
{\em regular}, and consequently the solutions  have
at worst logarithmic singularities, if and only if 
condition (\ref{regularity}) is satisfied.  Classically
this is called the {\em nonconfluent} case. This terminology 
is extended to the multivariate situation.
We shall assume throughout that the configuration $A$ is nonconfluent,
which implies that the $A$-hypergeometric system is
{\em regular holonomic} (\cite{gkz89},
\cite[2.4.11]{sst2}), and that it satisfies 
Assumptions~\ref{pyramid} and \ref{distinctpoints}. 
Since we are interested in the study of  rational
hypergeometric functions, we will also assume  that $\alpha \in
\Z^d$.

We recall that a configuration $A$ is said to be
Cayley if
$d =   2 r + 1$ and there exist vector configurations
$A_1,\dots,A_{r+1}$ in $\Z^r$ such that
\begin{equation}
\label{Cayley}
A \,\,\, =   \,\, \,
\{e_1\} \! \times A_1 \! \, \, \cup \,\,
\cdots \, \, \cup \,\,
\{e_{r+1}\} \! \times \! A_{r+1}
\,\, \subset \,\, \Z^{r+1} \times \Z^r ,
\end{equation}
where $e_1,\dots, e_{r+1}$ is the standard basis of
$\Z^{r+1}$.
Moreover, $A$ is said to be {\it essential} if the
Minkowski sum
$\,\sum_{i \in I} A_i \,$ has affine dimension at
least
$|I|$ for every proper subset $I$ of $\{1,\dots,r\}$. 

If $A$ is a codimension-two Cayley configuration then
the total number of points $n =  2r + 3$ and, if $A$ is
essential, each of the subsets $A_i$ must contain at least 
two points.
Hence, in an essential Cayley configuration of
codimension two, all but one of the
 $A_i$'s   contains two points
and the remaining one contains
three points.  A generic sparse
polynomial $f$ with support $A$ decomposes as
$$f(s_1, \dots, s_{r+1}, t_1, \dots, t_r) =   s_1
f_1(t) + \ldots +
s_{r+1}  f_{r+1}(t),$$
where $f_1, \dots, f_{r}$ are binomials with
respective supports
$A_1, \dots, A_r$ and $f_{r+1}$ is a trinomial with
support $A_{r+1}$.
Then, it is easy to deduce from \cite[\S5]{elim} that $A$
will be
essential if and only if $D_A \not=  1$, and
in this case $D_A$ equals the sparse resultant
$R_{A_1, \dots, A_{r+1}} (f_1, \dots, f_{r+1})$
(\cite{gkzbook}).

Conjecture~1.3 in \cite{rhf} asserts that the
only gkz-rational configurations are those affinely
isomorphic
to an essential Cayley configuration.  The
codimension-one case
as well as the dimension one, two, and three cases
have been
studied in \cite{rhf}.  In this section we verify this
conjecture for  the first significant open
case: codimension-two configurations in dimension four, i.e.
$n=7$, $d=5$.
We prove, moreover, that
for any such configuration the number of linearly
independent stable rational hypergeometric functions (see
Definition \ref{stable} below)
is $1$. In fact, a suitable derivative of any 
stable rational hypergeometric
functions  is a constant multiple of an explicit 
toric residue in the sense of
\cite{ compo, cox} associated to $A$ and the parameter
vector $\alpha \in \Z^{5}.$

\smallskip

The following result generalizes, in the case of
codimension-two 
configurations, Lemma~3.4 in \cite{rhf}.

\begin{lemma}\label{discspec} Let 
$A(I)$ be a non-splitting circuit of a
codimension-two configuration $A$.  Then there exists
$j\in I^c$ such that a positive power of
the discriminant $D_{A(I)}$ divides the
specialization $D_A\big|_{x_j=0}$.
\end{lemma}

\begin{proof}
We may assume without loss of generality that 
$I = \{1,\dots,r\}$ and that for $j>r$, $b_j =
\lambda_j e_2$,
$\lambda_j \in \Z$,  $\sum_{j>r} \lambda_j = \lambda >
0$.
This implies that $C(I^c) = \{b_{11},\dots,b_{r1}\}$ and
hence   $D_{A(I)}$ is, up to an integer constant, the
binomial
$$
\prod_{b_{i1}>0} b_{i1}^{b_{i1}}\cdot
\prod_{b_{i1}<0} x_{i}^{-b_{i1}} - 
\prod_{b_{i1}<0} b_{i1}^{-b_{i1}}\cdot
\prod_{b_{i1}>0} x_{i}^{b_{i1}},
$$
where the products run over $i\in I$.
On the other hand, it follows from \cite[\S 4]{elim}
that the discriminant $D(A)$ may be computed,
up to a monomial and integral factor, as the
resultant (with respect to $t$) of two polynomials
$p_1(t;x)$, $p_2(t;x)$ with the following
properties: $p_1(t;x)$ does not involve the
variables $x_j$, for $j>r$,  $p_1(0;x)= D_{A(I)}$,
and 
$$p_2(t;x) = a_1(t) t^\lambda m_1(x) - a_2(t) m_2(x), 
$$
with $$m_1(x) = \prod_{b_{k2}<0} x_{k}^{-b_{k2}},\ 
m_2(x) = \prod_{b_{k2}>0} x_{k}^{b_{k2}}, $$
where the products run over $k=1,\dots,n$.

By assumption, there exists $j>r$ such that $b_{j2} =
\lambda_j > 0$.
Hence, setting $x_j = 0$ we obtain:
\begin{eqnarray}\label{discformula}
D_A\big|_{x_j=0}& = & {\rm
Res}_t\bigl(p_1(t;x)\big|_{x_j=0},
p_2(t;x)\big|_{x_j=0}\bigr)\\
& =&  {\rm Res}_t (p_1(t;x),a_1(t) t^\lambda m_1(x)).
\nonumber
\end{eqnarray}
Poisson's formula for resultants now implies that the
resultant
$$ {\rm Res}_t (p_1(t;x), t^\lambda)\  =\ 
p_1(0;x)^\lambda
\ =\ D_{A(I)}^\lambda  $$
divides (\ref{discformula}).
\end{proof}

\bigskip

\begin{theorem}\label{generic}
A codimension-two configuration which contains an
unbalanced non-splitting circuit is not gkz-rational.
\end{theorem}

\begin{proof}
This result generalizes \cite[Theorem~1.2]{rhf} whose
proof we follow.
Let $A =   \{a_1,\dots,a_n\}$ be a gkz-rational
configuration of 
codimension two.  Suppose
$f =   P/Q$ is  a rational $A$-hypergeometric function
of degree $\alpha\in\Z^d$,
where $P, Q \in \C[x_1,\dots,x_n]$ are relatively
prime. Assume, moreover, that
 $D_A$ is not a monomial and divides $Q$.

We claim that any non-splitting circuit $A(I)$ of $A$
is balanced.
We may assume that $I = \{1,\dots,r\}$, $r<n$ and,
because of Lemma~\ref{discspec}, that 
$D_A\big|_{x_n=0}$  is not a monomial.
Set $t=  x_n$,
$\tilde A =   \{a_1,\dots,a_{n-1}\}$, $\tilde x =
(x_1,\dots,x_{n-1})$. We expand the 
$A$-hypergeometric function
$f$ as
\begin{equation}
\label{SeriesForMircea}
f(\tilde x;t) \quad =   \quad
\sum_{\ell\geq \ell_0}\,R_\ell(\tilde x)\cdot
t^\ell\,,
\end{equation}
where each $R_\ell(\tilde x)$ is a rational $\tilde
A$-hypergeometric
function of degree $\alpha - \ell\cdot a_n$.
Since  $D_A\big|_{x_n=0}$  is not a monomial, 
it follows from \cite[Lemma~3.3]{rhf} that some 
coefficient $R_\ell(\tilde x)$
is not a Laurent polynomial. Hence, $\tilde A$ is
gkz-rational.
But $\tilde A$ is a codimension-one configuration;
indeed, it is a
pyramid over the circuit $A(I)$. Thus, it follows from
\cite[Theorem~2.3]{rhf} that $A(I)$ must be balanced.
\end{proof}

\medskip

\begin{theorem}\label{th:conj}
A codimension two, dimension four, gkz-rational
configuration is affinely equivalent to
an essential Cayley configuration.
\end{theorem}

\begin{proof} 
A Gale dual of $A$ is a configuration $B$ of seven
lattice vectors in the plane. Hence, we may apply
the classification given in Theorem~\ref{classification}.
If $B$ is as in a), it contains an unbalanced non-splitting
circuit and then
it may not be dual to a gkz-rational configuration by
Theorem~\ref{generic}.
A direct computation
(see also \cite[Corollary 4.5]{elim}) shows that the
discriminant of any
configuration $A$
whose Gale dual satisfies c) must be $1$ and hence 
$A$ is not gkz-rational. Finally, note that no lattice
configuration may be as in b) or e).

We are left with Gale dual configurations $B =   B_1
\cup B_2 \cup B_3$ as in d).
Two of the subconfigurations contain $2$ vectors and
the third has $3$ vectors not on a line. Then,  $A$ is
Cayley essential and by
\cite[Theorem~1.5]{rhf} it is gkz-rational.
\end{proof}

\medskip
We will now construct rational hypergeometric
functions associated
with an essential Cayley configuration.
Let $A$ be  Cayley essential  with $r=  2, \, n=  7$.
We may assume,
modulo a change of coordinates, that
the columns of $A$ index the coefficients in
a sparse system of $3$ equations of the form:
\begin{eqnarray} \label{cayley7}
f_1 \quad =   & x_1 \, + \, y_1\cdot t_1^{\gamma_1}
\qquad\nonumber\\
f_2 \quad =   & x_2 \, + \, y_2 \cdot t_2^{\gamma_2}
\qquad \\
f_3 \quad =   & x_3 \, + \,
y_3 \cdot  t_1^{\alpha_1} t_2^{\alpha_2} \,+ \,
z_3 \cdot  t_1^{\beta_1}  t_2^{\beta_2}  ,\nonumber
\end{eqnarray}
with $\gamma_1, \gamma_2 \in \Z_{>0}$ and
$(\alpha_1, \beta_1) , (\alpha_2, \beta_2)$  non-zero
integral vectors.
For $i\not= j$, let $V_{ij}$ denote the common zero
set of $f_i$ and $f_j$
in the torus $(\C^*)^2$.
For any $a \in \Z^2, c\in \Z_{>0}^3,$ the sum
$R_{ij}(c,a)$ of the
local  residues
at all the points in $V_{ij}$ of the differential
$2$-form
$$  \frac{t^{a} }
{ f_1^{c_1} f_2^{c_2} f_3^{c_3}} \frac{dt_1}{t_1}
\wedge
\frac{dt_2}{t_2}, $$
is a rational $A$-hypergeometric function of degree
$(-c, -a)$.
  When $a$ lies in the interior of the Minkowski
sum of the Newton polygons of $f_1^{c_1},f_2^{c_2}$
and $f_3^{c_3}$,
 we have by Theorem~4.12 in \cite{compo} that
\begin{equation}\label{stokes}
R_{12}(c,a) =   - R_{13}(c,a) =   R_{23}(c,a).
\end{equation}
Call  this rational function $R(c,a)$.  It has an
integral representation
(\cite{compo,tsikh}):
$$ R(c,a)\ =\ \int_\Gamma \ \frac{t^{a} }
{ f_1^{c_1} f_2^{c_2} f_3^{c_3}} \frac{dt_1}{t_1}
\wedge
\frac{dt_2}{t_2}$$
for an appropriate real $2$-cycle $\Gamma$ in the torus.

As shown in \cite[Theorem~7]{mega96}, we can
differentiate $R(c,a)$ under the integral sign,
for any $c\in \Z_{>0}^3, a \in \Z^2$.  We have, for
example,
\begin{equation}\label{derx}
\partial_{x_1} R(c,a)\  =   \
 - c_1 \, R(c+ (1,0,0),a), \hbox{\ and}
\end{equation}
\begin{equation}\label{dery}
\partial_{z_3} R(c,a)  =   \  - c_3 \,
R(c + (0,0,1), a + \beta).
\end{equation}

\begin{definition} \label{stable}
A rational function $f$ is called {\it unstable} if it
is annihilated
by some iterated derivative.
Otherwise we say that $f$ is {\it stable}.
\end{definition}
Thus, a rational function $f$ is unstable if it is a
linear combination of rational functions
that depend  polynomially on at least one of the
variables.

\begin{proposition}  \label{prop:nonzero}
Let $A$ be a Cayley essential configuration with $r= 
2,\, n=  7$ and let  $c \in \Z^3_{>0}$.
Assume the integer vector $(c,a)$ lies in the integer
image of $A$. Then, $R(c,a) \not=   0.$ Moreover,
$R(c,a)$ is stable.
\end{proposition}

\begin{proof}
We argue as in the proof of Proposition~4.4 in
\cite{binom} which asserts a similar statement in case
$f_1,f_2,f_3$ are binomials.  

By definition, $R(c,a)$ may be computed as
the sum of the local residues over $V_{12}$ of the
form
$$  \frac {(t^{a}/f_3^{c_3})}
{ f_1^{c_1} f_2^{c_2}} \frac{dt_1}{t_1} \wedge
\frac{dt_2}{t_2}$$
with respect to $(f_1, f_2)$.
By the derivative formulae such as (\ref{derx}), it
suffices to show the result for
$c_1=  c_2=  c_3 =  1$.
We may expand the numerator $t^a/f_3$ as a Laurent
series of the form
\begin{equation*}
\frac {t^{a} }  {f_3} =    \frac 1 {x_3} \, \frac
{t^a} { 1 + (y_3/x_3)
t^\alpha + (z_3/x_3) t^\beta} =
\sum_{m\in \N^2} c_m \frac {y_3^{m_1}
z_3^{m_2}}{x_3^{m_1+m_2+1}}
t^{a+m_1 \alpha + m_2 \beta},
\end{equation*}
whith $c_m \ne 0$ for all $m =   (m_1,m_2) \in \N^2$. 
 By Lemma~4.2 in
\cite{binom}, the global residue  with respect
to $f_1,f_2$ of the form
$$ \frac {t^{a + m_1 \alpha + m_2 \beta}} {f_1 \cdot
f_2} \
\frac{dt_1}{t_1} \ \wedge
\ \frac {dt_2}{t_2}$$
vanishes unless there exists
$\nu(m)=(\nu_1(m),\nu_2(m)) \in\Z^2$
such that
$$a + m_1 \alpha + m_2 \beta  \ =\  (\nu_1(m)
\gamma_1, \nu_2(m)
\gamma_2).$$
In this case the residue
is equal (up to sign) to the Laurent monomial
$$x_1^{\nu_1(m)-1} y_1^{-\nu_1(m)}
x_2^{\nu_2(m)-1} y_2^{-\nu_2(m)}. $$
 Exchanging the residue integral and the sum of the
series,  we have
\begin{equation*}
R(c,a) =   \sum \pm
 c_m \frac {y_3^{m_1}
z_3^{m_2}}{x_3^{m_1+m_2+1}}x_1^{\nu_1(m)-1}
y_1^{-\nu_1(m)}
x_2^{\nu_2(m)-1} y_2^{-\nu_2(m)},
\end{equation*}
where the sum runs over all $m\in \N^2$ for which the
equation
\begin{equation}\label{resequation}
{a + m_1 \alpha + m_2 \beta = (\nu_1(m) \gamma_1,
\nu_2(m)
\gamma_2)}
\end{equation}
has an integral solution $\nu(m)$.

Since $(c,a)$ lies in the integer image of $A$,
(\ref{resequation}) has a solution for one,
and a fortiori infinitely many pairs $(m_1,m_2)$.
Moreover, the assumption that $\alpha$ and $\beta$
cannot be both multiples of $e_1$ or $e_2$
ensures that there will be infinitely many values
of $\nu_1(m)$ and $\nu_2(m)$ as well. Hence,
$R(c,a)$ and all its iterated derivatives are non zero.
\end{proof}


\begin{remark} \label{EJ}
We recall that  the open cone in $\R^5$,
$${\mathcal E}
\  =   \  \left\{\sum_{i=1}^7 \nu_i a_i\,:\,\nu_i\in
\R,\,\nu_i<0\right\}
$$
is called
the {\it Euler-Jacobi \/} cone of $A$.
Then $(-c,-a)$ lies in ${\mathcal E}$ if and only if
$c_1,c_2,c_3 >0$ and
$a$ lies in the interior of the Minkowski sum of the
convex hulls of $c_1 A_1, c_2 A_2$ and $c_3 A_3$ (i.e.
the Newton polygons of $f_1^{c_1}, f_2^{c_2}, f_3^{c_3}$).
Note that given  $\alpha$ in the integer image
of $A$, the vector $\alpha - A \cdot  v$ lies in 
${\mathcal E}$ for any
integer vector $v \in \N^n$ with $v_i \gg 0$ for all
$i=  1,\dots,n.$
\end{remark}

\begin{theorem} \label{th:dim1}
Let $A \subset \Z^{5\times 7}$ be a Cayley essential
configuration. Let $(-c,-a) \in \EE \cap A \cdot \Z^7$.
Then, the  dimension of the space of rational
$A$-hypergeometric functions of degree $(-c,-a)$ is
equal to $1$ and it is spanned by $R(c,a)$.
\end{theorem}

\begin{proof}
By Proposition \ref{prop:nonzero}, it is enough to
show that the
dimension of the space of rational $A$-hypergeometric
functions cannot exceed $1$.
In order to prove this statement we need to recall the
characterization
of logarithm-free $A$-hypergeometric series of degree
$(-c,-a)$.  We refer to
\cite[\S3.4]{sst2} for details and proofs.

Let $B \in \Z^{7\times 2}$ be a Gale dual of $A$; as
before, we also denote by
$B =  \{b_1,\ldots,b_7\}$ the planar configuration of
its row-vectors. Set
$\LL = \ker_\Z(A)$ and
$\LL_\R = \LL\otimes_\Z \R$. 
For $v \in \Z^7$  such that $A \cdot v =   (-c,-a)$,
the plane
$v + \LL_\R$ may be identified with $\R^2$ via the
affine isomorphism
\begin{equation}\label{isom}
\lambda \in \R^2 \ \mapsto \ v \ + \ B \cdot \lambda \,.
\end{equation}
This correspondence maps $\Z^2$ to  $v + \LL$.
 Consider the affine
arrangement of hyperplanes $\mathcal H$ in $\R^2$
given by
\begin{equation}\label{hyperpl}
 H_j :=  \{\lambda \in \R^2\ :\ \langle b_j,\lambda
\rangle
=   -v_j\},\ \  j=  1,\dots,7.
\end{equation}
Each $H_j$ is oriented by the choice of normal vector
$b_j$.
The negative support of a vector $\lambda$ is defined
as the set of all the indices $j =  1,\dots,7 \,$ for
which
$\langle b_j,\lambda \rangle < -v_j$. The (convex) set
of points
with the same negative support is called a cell of
the hyperplane arrangement.  A cell $\Sigma$ is
minimal
if $\Sigma \cap \Z^2 \ne \emptyset$ and the negative
support of the elements in this set is minimal with
respect
to inclusion among the supports of integer points in
any
other cell. If $(-c,-a)\in {\mathcal E} $ then all minimal cells
are unbounded
and correspond to minimal cells in the oriented
central arrangement.
To each minimal cell we may attach a formal
Laurent series $\phi_\Sigma$, unique up to constant, 
whose exponents are the points
in $\Sigma \cap \Z^2$, and which is  annihilated by
the action of the
$A$-hypergeometric ideal $H_A((-c,-a))$.

Recall that we may assume that $B$ satisfies
$$ b_1+b_2\ =\ b_3+b_4\ =\ b_5+b_6 +b_7 \ =\ 0.$$
The central arrangement defined by $b_5, b_6, b_7$ has
three minimal cells.  The lines orthogonal to $b_1$, $b_3$, may
intersect at most
two of those of cells, hence it is clear that there
exists $w \in \R¾$ not parallel to any vector in $B$, for
which there is a single minimal 
cell  $\Sigma$ lying in a half-space
$\langle w,\lambda \rangle > \rho,$ for some $\rho \in \R$.

To complete the proof we need to show that the
rational function
$R(c,a)$ has a Laurent expansion 
with support in $\Sigma$.  Moreover, since in the
Euler-Jacobi cone there are no bounded minimal regions
(cf. \cite[Proposition~2.6]{binom}) it is enough to show
that the linear functional $\langle w,.\rangle$ is bounded
below in the preimage under (\ref{isom}) of the
support of a Laurent expansion of $R(c,a)$.

Write $R(c,a) =   P/Q$ with $P,Q$
 polynomials without common factors.
Since the discriminant $D_A$ divides $Q$, the Newton
polytope of $Q$,
$N(Q)$, must be two-dimensional.
Moreover, since $Q$ is $A$-homogeneous,
$N(Q)$ lies in a plane $\Pi$ parallel to $\LL_\R$.
Let $\tilde w\in \LL_{\R}$ be the
unique element such that $B\cdot \tilde w = w$ and
choose $\mu_{0}\in \Pi$.
Perturbing $w$ if necessary, we may assume that the
linear functional
$$\mu\in \Pi \mapsto \langle \tilde w , (\mu -
\mu_{0})\rangle$$
acting on $N(Q)$ attains its minimum only at a vertex
$\nu_0$ of $N(Q)$.  This is
clearly independent of the choice of $\mu_{0}$, hence
we may assume
from now on that $\mu_{0 } = \nu_{0}$.
Let $\nu_1,\dots,\nu_a$ denote the remaining integral
points of $N(Q)$; note that, by construction,
$$\langle \tilde w,(\nu_{j}-\nu_{0}) \rangle\ >  0\,;\
\ j=1,\dots,a.$$
If we expand $1/Q$ at $\nu_0$ (cf. \cite[6.1.b]{gkzbook})
we obtain
$$ 1/Q \ =\ x^{-\nu_0}\cdot\sum_{\alpha\in \N^a}
c_\alpha
x^{\sum\alpha_j(\nu_j - \nu_0)}.$$
Now, for any monomial $x^u$ in $P$, the exponents of
$x^u/Q$ are of the
form $(u-\nu_0) + \sum\alpha_j(\nu_j -
\nu_0)$ and clearly, $\langle w,.\rangle$ is bounded
below in in the preimage under (\ref{isom}) of this
set.  Finally, since the support of $P/Q$ is contained
is a finite union of sets of this form, the result follows.
\end{proof}

As a Corollary, we identify all stable rational
functions for any
integer homogeneity, proving in this case 
Conjecture~5.7  in \cite{rhf}.

\begin{corollary}
Let $\alpha\in\Z^5$. The dimension of
 the space of rational $A$-hypergeometric functions
of degree $\alpha,$
modulo  unstable ones, is at most  $1$. Moreover,  for
any rational $A$-hypergeometric  function,
an iterated derivative of it is zero or a multiple of
a toric residue.
\end{corollary}

\begin{proof}
If there is a non-zero rational $A$-hypergeometric function
$f$ of degree $\alpha$, then $\alpha$ is necessarily
in the integer image of $A$. Moreover,
there exists $v\in \N^7$ such
that $\alpha - A\cdot v \in \EE \cap A \cdot \Z^7$.
Hence, the derivative $D_v(f)$ 
is a multiple (possibly zero) of a toric residue.
\end{proof}

\smallskip

\begin{example}  We recall that there are $14$
complete hypergeo\-me\-tric
Horn series in two variables and of {\em order} two
\cite[\S 5.7.1]{erdelyi}.  Modulo monomial changes of
variables and analytic continuation, they give rise to
$8$ different
functions \cite[\S 3.2]{gkz89}.  In combinatorial
terms they correspond to  Gale configurations in $\Z^2$ of
vectors with the property that the sum of positive first or
second entries is  equal to two.  In particular,
there exists a unique such configuration of seven
vectors and it is associated with
the classical hypergeometric series $F_2$, $F_3$ and
$H_2$.
The configuration $A$ and a Gale dual $B$  can be 
chosen as:
\begin{equation*}
A \  =   \
\left(
\begin{array}{rrrrrrr}
1 & 1 & 0 & 0 & 0 & 0 & 0 \\
0 & 0 & 1 & 1 & 0 & 0 & 0 \\
0 & 0 & 0 & 0 & 1 & 1 & 1 \\
0 & 1 & 0 & 0 & 0 & 1 & 0 \\
0 & 0 & 0 & 1 & 0 & 0 & 1
\end{array}
\right) \, \  ;\quad
B \  =   \
\left(
\begin{array}{rr}
\ 1 & \ 0  \\
-1 & \ 0  \\
\ 0 & \ 1  \\
\ 0 & -1  \\
-1 & -1 \\
\ 1 & \ 0  \\
\ 0 & \ 1
\end{array}
\right)\,.
\end{equation*}
Thus, this is a Cayley essential configuration, with
$A_1 = \{0,e_1\},\, A_2 = \{0,e_2\}, \, A_3 =
\{0,e_1,e_2\}\subset \Z^2$.
Let $f_1(t), f_2(t), f_3(t)$ be as in (\ref{cayley7}),
$c=(1,1,1)$, and
$a=(0,0)$.  Then $R_{12}(c,a)(x_1,\dots,z_3)$ may be computed as a single
local residue
\begin{eqnarray*}
R_{12}(c,a) & = &
{\rm
Res}_{(-\frac{x_1}{y_1},-\frac{x_2}{y_2})}\left(\frac
{1/(t_1
t_2 (x_3 + y_3 t_1 + z_3 t_2))}{(x_1 + y_1 t_1)(x_2 +
y_2 t_2)} \,
dt_1\wedge
dt_2\right)\\
&=\ & \frac {y_1y_2}{x_1x_2(y_1y_2x_3 -x_1y_2y_3 -
y_1x_2z_3)}\,.
\end{eqnarray*}
Although (\ref{stokes}) does not hold in this degree,
this is the unique, up to constant, stable
rational $A$-hypergeo\-me\-tric
function of degree $(-c,-a)$.  The function
$R_{12}(c,a)$ may
also be written as
$$R_{12}(c,a) = \frac {1}{x_1x_2x_3}\,\left(
\frac {1}{1-u-v}
\right)\,;
\ u = \frac {x_1 y_3}{y_1 x_3}\,, \ v = \frac {x_2
z_3}{y_2 x_3}$$
and
$$\frac {1}{1-u-v} = \sum_{m,n\geq 0} \binom{m+n}{m}
u^m v^n =
F_2(1,\beta,\beta',\beta,\beta';u,v)\, ,$$
for any $\beta, \beta' \in \C \, \setminus \, \Z_{\leq
0}.$
Note that the vector $(-1,-1,-1,0,0)$ does not lie in the
Euler-Jacobi cone and, indeed, there are
unstable rational
functions of this degree. They are also realized as
residues:
\begin{eqnarray*}
R_1(x_1,\dots,z_3) & = &
 {\rm Res}_{(-\frac{x_1}{y_1},0)}\left(\frac {1/(t_1
f_2(t) f_3(t))}{t_2(x_1 + y_1 t_1)} \, dt_1\wedge
dt_2\right)\\
&=\ & \frac {y_1}{x_1x_2(y_1x_3 -x_1y_3)}\,,
\end{eqnarray*}
\begin{eqnarray*}
R_2(x_1,\dots,z_3)  & = &
 {\rm Res}_{(0,-\frac{x_2}{y_2})}\left(\frac {1/(t_2
f_1(t) f_3(t))}{(x_2 + y_2 t_2)t_1} \, dt_1\wedge
dt_2\right)\\
&=\ & \frac {y_2}{x_1x_2(y_2x_3 -x_2z_3)}\,, \hbox{\
and}
\end{eqnarray*}
\begin{eqnarray*}
R_3(x_1,\dots,z_3) & = &
 {\rm Res}_{(0,0)}\left(\frac {1/(f_1(t)
f_2(t) f_3(t))}{t_1t_2} \, dt_1\wedge
dt_2\right) \\
& = &\ \frac 1{x_1x_2x_3}\,.
\end{eqnarray*}
These four rational functions are linearly
independent and  the
holonomic rank of the hypergeometric system
$H_A(-c,-a)$ is
$4$, so in this example all solutions are rational,
which is highly exceptional.

On the other hand, for  $c = (1,1,2)$, $a=(1,1)$, the
vector $(-c,-a)$ is in $\EE$ and there
exists a unique, up to constant, rational
$A$-hypergeo\-me\-tric function of this degree:
\begin{eqnarray*}
R(c,a)  & = &
 {\rm Res}_{(-\frac{x_1}{y_1},-\frac{x_2}{y_2})}\left(\frac
{1/
(x_3 + y_3 t_1 + z_3 t_2)^2}{(x_1 + y_1 t_1)(x_2 + y_2
t_2)} \, dt_1\wedge
dt_2\right)\\
&=\ & \frac {y_1y_2}{(y_1y_2x_3 -x_1y_2y_3 -
y_1x_2z_3)^2}\\
&=\ & \frac {1}{x_3^2 y_1 y_2} \cdot
 \frac1{(1-u-v)^2} \\
& = &\  \frac {1}{x_3^2 y_1 y_2} \cdot
F_2(2,1,\beta',2,\beta';u,v) \,,
\end{eqnarray*}
for any parameter $\beta' \in  \C \,\setminus \,
\Z_{\leq 0}.$
\end{example}


\begin{thebibliography}{99}

\bibitem{matroids} A.~Bj\"orner, M.~Las~Vergnas,
B.~Sturmfels, N.~White, and G. Ziegler:
Oriented Matroids, Cambridge University Press, 1993.

\bibitem{compo} E.~Cattani, D.~Cox, and A.~
Dickenstein: {Residues in
toric varieties}. {\it Compositio Mathematica~}
{\bf 108} (1997), 35--76.

\bibitem{cdd} E.~Cattani, C.~D'Andrea, and
A.~Dickenstein:
{The $\mathcal A$-hypergeometric system associated
with a monomial curve.} {\it Duke Math.~J.~} {\bf 99}
(1999), 179--207.


\bibitem{mega96}
E.~Cattani and A.~Dickenstein:
{A global view of residues in the torus}.
{\it Journal of Pure and Applied Algebra~} {\bf 117 \&
118} (1997), 119--144.

\bibitem{rhf} E.~Cattani, A.~Dickenstein, and
B.~Sturmfels: {Rational hypergeometric functions}.
{\it Compositio Mathematica~}{\bf 128} (2001),
217--240.

\bibitem{binom} E.~Cattani, A.~Dickenstein, and
B.~Sturmfels: {Binomial Residues}. {\it Ann. Inst. Fourier~}
{\bf 52} (2002), 687--708. 

\bibitem{cox} D.~Cox:
{ Toric residues}. {\it Arkiv f\"or Matematik~} {\bf
34} (1996), 73--96.

\bibitem{elim} A.~Dickenstein and
B.~Sturmfels: {Elimination theory in codimension two}.
{\it Journal of Symbolic Computation} {\bf 34} (2002), 119--135.

\bibitem{erdelyi} A. Erd{\'e}lyi et al. :
{\it Higher Transcendental Functions.}
{Based, in part, on notes left by Harry Bateman},
 {McGraw-Hill Book Company, Inc., New
York-Toronto-London}, {1955}.


\bibitem{gkz89} I. M.~Gel'fand, A.~Zelevinsky, and
M.~Kapranov:
{Hypergeometric functions and toral manifolds},
{\it Functional Analysis and its Appl.}
{\bf 23} (1989), 94--106.

\bibitem{gkz90} I. M.~Gel'fand, M.~Kapranov, and
A.~Zelevinsky:
{Generalized Euler integrals and ${\mathcal
A}$-hypergeometric functions}, {\it Advances in Math.} {\bf 84} 
(1990), 255--271.

\bibitem{gkzbook} I. M.~Gel'fand, M.~Kapranov, and
A.~Zelevinsky:
{\it Discriminants, Resultants and Multidimensional
Determinants,} Birkh\"auser, Boston, 1994.

\bibitem{hille} E.~Hille:
{\it Ordinary Differential Equations in the Complex
Domain,} Dover, N.Y., 1997.

\bibitem{ressayre} N.~Ressayre: {Balanced configurations
of $2n+1$ plane vectors}. Preprint, (2002),
arXiv:math.RA/0206234.


\bibitem{sst}
M.~Saito, B.~Sturmfels, and N.~Takayama:
{Hypergeometric polynomials and integer programming}.
{\it Compositio Mathematica~}{\bf 115} (1999), 185--204.

\bibitem{sst2}
M.~Saito, B.~Sturmfels, and N.~Takayama:
{\it Gr\"obner Deformations of Hypergeometric
Differential Equations},
Algorithms and Computation in Mathematics, Volume {\bf
6}, Springer-Verlag, Heidelberg, 1999.

\bibitem{tsikh} A.~Tsikh:
{\it Multidimensional Residues and Their
Applications}, American Math. Society, Providence, 1992.

\end{thebibliography}
\end{document}